\documentclass{article}

\usepackage{amssymb,amsmath,color, xcolor}
\usepackage{amsthm,setspace}
\usepackage{url}
\usepackage{tikz,subcaption}
\usepackage[enableskew]{youngtab}
\usepackage{ytableau,varwidth}
\definecolor{dartmouthgreen}{rgb}{0.05, 0.5, 0.06}

\definecolor{red}{rgb}{1,0,0}
\definecolor{blue}{rgb}{.2,.2,.8}

\newtheorem{theorem}{Theorem}[section]
\newtheorem{corollary}[theorem]{Corollary}

\theoremstyle{definition}
\newtheorem{definition}{Definition}

\newtheorem{remark}{Remark}

\newcommand*\samethanks[1][\value{footnote}]{\footnotemark[#1]}

\begin{document}

	\title{Dyson's crank and unimodal compositions}
\author{Cristina Ballantine
	\\
	\footnotesize Department of Mathematics and Computer Science, College of The Holy Cross\\
	\footnotesize Worcester, MA 01610, USA \\
	\footnotesize cballant@holycross.edu
	\and Mircea Merca\samethanks
	\\ 
	\footnotesize Department of Mathematics, University of Craiova, DJ 200585 Craiova, Romania\\
	\footnotesize mircea.merca@profinfo.edu.ro
}
	\date{}
	\maketitle


\begin{abstract}
The crank is a partition statistic requested by Dyson in 1944 in order to combinatorially prove a Ramanujan congruence of Euler's partition function $p(n)$. In this paper, we 
provide connections between Dyson's crank and unimodal compositions. Somewhat unrelated, we give a combinatorial proof of a new truncated Euler pentagonal number theorem due to Xia and Zhao. 
\\
\\
{\bf Keywords:} crank, integer partitions, unimodal compositions
\\
\\
{\bf MSC 2010:}  11P81, 11P82, 05A19, 05A20 
\end{abstract}

\section{Introduction}

Recall that a composition of a positive integer $n$ is a sequence $\lambda=(\lambda_1,\lambda_2,\ldots,\lambda_k)$ of positive integers whose sum is $n$, i.e.,
\begin{equation}\label{EQ1}
n=\lambda_1+\lambda_2+\cdots+\lambda_k.
\end{equation}
The positive integers in the sequence are called parts \cite{Andrews98}. 
The composition $(\lambda_1,\lambda_2,\ldots,\lambda_k)$ is unimodal if there exist an integer $t\in\{1,2,\ldots,k\}$ such that
$$ \lambda_1 \leqslant \lambda_2 \leqslant \cdots \leqslant \lambda_t \quad\text{and}\quad 
\lambda_t \geqslant \lambda_{t+1} \geqslant \cdots \geqslant \lambda_{k}.$$

\begin{definition}
	For positive integers $n$ and $m$,  we define:
	\begin{enumerate}
		\item [(i)] $u_0(n)$ to be the number of unimodal compositions of $n$. We set $u_0(0)=0$. 
		\item [(ii)] $u_m(n)$ to be the number of unimodal compositions of $n+m$ where the maximal part appears exactly $m$ times. Note that $u_m(0)=1$.
	\end{enumerate}	 
\end{definition}

For example, there are $u_0(4)=8$ unimodal compositions of $4$:
\begin{align*}
& (1,1,1,1),\ (1,1,2),\ (1,2,1),\ (1,3),\ (2,1,1),\ (2,2),\ (3,1),\ (4).
\end{align*}
There are $u_1(4)=12$ unimodal compositions of $5$ where the maximal part appears once:
\begin{align*}
& (1,1,1,2),\ (1,1,2,1),\ (1,1,3),\ (1,2,1,1),\ (1,3,1)\\
& (1,4),\ (2,1,1,1),\ (2,3),\ (3,1,1),\ (3,2),\ (4,1),\ (5).
\end{align*}
There are $u_2(4)=4$ unimodal compositions of $6$ where the maximal part appears exactly twice:
\begin{align*}
& (1,1,2,2),\ (1,2,2,1),\ (2,2,1,1),\ (3,3).
\end{align*}

In this article we prove the following identity both analytically and combinatorially.

\begin{theorem}\label{Th1}
	For $n\geqslant 0$, $u_0(n)=u_1(n)-u_2(n)$.
\end{theorem}

When the order of the integers $\lambda_i$ does not matter, the representation \eqref{EQ1} is known as an integer partition. For consistency, a partition of $n$ will be written with the summands in non-increasing order \cite{Andrews98}. For example, the following are the partitions of $5$:
\begin{equation}
(5),\ (4,1),\ (3,2),\ (3,1,1),\ (2,2,1),\ (2,1,1,1),\ (1,1,1,1,1). \label{p5}
\end{equation}
If $\lambda$ is a composition or a partition whose parts sum up to $n$, we say that the  size of $\lambda$ is $n$ and write $|\lambda|=n$. We denote the number of parts of $\lambda$ by $\ell(\lambda)$.
In 1988, Andrews and Garvan \cite{Andrews88} defined the crank of an integer partition as follows. 
The crank of a partition is the largest part of the partition if there are no ones as parts and otherwise is the number of parts larger than the number of ones minus the number of ones. More precisely,
for a partition $\lambda=(\lambda_1,\lambda_2,\ldots,\lambda_k)$,
 let 
$\omega(\lambda)$ denote the number of $1$'s in $\lambda$, and $\mu(\lambda)$ denote the number of parts of $\lambda$ larger than $\omega(\lambda)$. Note that $\lambda_1$ is the largest part of $\lambda$. Then, the  crank, $c(\lambda)$, of $\lambda$ is defined  by
$$
c(\lambda) = \begin{cases}
\lambda_1, & \text{if $\omega(\lambda)=0,$}\\
\mu(\lambda)-\omega(\lambda),& \text{if $\omega(\lambda)>0.$}
\end{cases}
$$

\begin{definition}
	For a non-negative integer $n\neq 1$,  we define: 
	\begin{enumerate}
		\item [(i)] $C_0(n)$ to be the number of partitions $\lambda$ of $n$ with $c(\lambda)> 0$.
		\item [(ii)] $C_1(n)$ to be the number of partitions $\lambda$ of $n$ with $c(\lambda)\geqslant 0$;
		\item [(iii)] $C_2(n)$ to be the number of partitions $\lambda$ of $n$ with $c(\lambda)= 0$.
	\end{enumerate}
	We set $C_0(1)=1, C_1(1)=0$, and $C_2(1)=-1$.
\end{definition}

It is clear that $C_2(n)=C_1(n)-C_0(n)$. 
For instance, there are $C_0(5)=3$ partitions of $5$ with positive crank:
\begin{align*}
& (5),\ (3,2),\ (2,2,1).
\end{align*}
There are $C_1(5)=4$ partitions of $5$ with non-negative crank:
\begin{align*}
& (5),\ (4,1),\ (3,2),\ (2,2,1).
\end{align*}

In \cite{Andrews12}, while investigating the truncated Euler's pentagonal number theorem, Andrews and Merca introduced the partition function $M_k(n)$, which counts the number of partitions of $n$ where $k$ is the least positive integer that is not a part and there are more parts $>k$ than there are parts $<k$.  For instance, we have $M_3(18)=3$ because the three partitions in question are 
$$(5,5,5,2,1),\ (6,5,4,2,1),\ (7,4,4,2,1).$$ 

In this paper, we prove the following result connecting the partition functions $u_m(n)$, $C_m(n)$ and $M_k(n)$ when $m\in\{0,1,2\}$.

\begin{theorem}\label{Th2}
	Let $k$ and $n$ be positive integers. For $m\in\{0,1,2\}$, we have
	\begin{align*}
	& (-1)^{k} \left(C_m(n)-\sum_{j=1-k}^k (-1)^j u_m\big(n-j(3j-1)/2\big) \right)
		= \sum_{j=0}^n C_m(j)\,M_k(n-j).
	\end{align*}
\end{theorem}

An immediate consequence of this theorem  is  the following infinite
family of linear inequalities.

\begin{corollary}
	Let $k$ and $n$ be positive integers. For $m\in\{0,1,2\}$, we have	
	\begin{align*}
	& (-1)^{k} \left( C_m(n) - \sum_{j=1-k}^k (-1)^j u_m\big(n-j(3j-1)/2\big)   \right) \geqslant 0,
	\end{align*}
	with strict inequality if $n\geqslant k(3k+1)/2$.
\end{corollary}

For example, some special cases of this corollary are:
\begin{align*}
& u_m(n)-u_m(n-1)\geqslant C_m(n),\\
& u_m(n)-u_m(n-1)-u_m(n-2)+u_m(n-5)\leqslant C_m(n),\\
& u_m(n)-u_m(n-1)-u_m(n-2)+u_m(n-5)+u_m(n-7)-u_m(n-12)\geqslant C_m(n).
\end{align*}

Very recently, Xia and Zhao \cite{Xia} defined $\widetilde{P}_k(n)$ to be the number of partitions of $n$ 
in which every part $\leqslant k$ appears at least once
and the first part larger that $k$ appears at least $k + 1$ times. For example, $\widetilde{P}_2(17)=9$, and the partitions in question are: 
\begin{align*}
& (5,3,3,3,2,1), (4,4,4,2,2,1), (4,4,4,2,1,1,1),\\
& (4,3,3,3,2,1,1), (3,3,3,3,2,2,1), (3,3,3,3,2,1,1,1),\\
& (3,3,3,2,2,2,1,1), (3,3,3,2,2,1,1,1,1), (3,3,3,2,1,1,1,1,1,1).
\end{align*}

In analogy with Theorem \ref{Th2}, we have the following result.

\begin{theorem}\label{Th3}
	Let $k$ and $n$ be positive integers. For $m\in\{0,1,2\}$, we have
	\begin{align*}
	& (-1)^{k-1}\left( C_m(n)-\sum_{j=-k}^k (-1)^j u_m\big(n-j(3j-1)/2\big) \right) 
		= \sum_{j=0}^n C_m(j)\,\widetilde{P}_k(n-j).
	\end{align*}	
\end{theorem}

An immediate consequence of this theorem  is  the following infinite
family of linear inequalities.

\begin{corollary}
	Let $k$ and $n$ be positive integers. For $m\in\{0,1,2\}$, we have	
	\begin{align*}
	& (-1)^{k-1} \left( C_m(n) - \sum_{j=-k}^k (-1)^j u_m\big(n-j(3j-1)/2\big)   \right) \geqslant 0.
	\end{align*}
\end{corollary}

When $k\to\infty$, Theorem \ref{Th2} or \ref{Th3} yields the following decomposition of $C_m(n)$ for $m\in\{0,1,2\}$.

\begin{corollary}\label{decomp}
	Let $n$ be a positive integer. For $m\in\{0,1,2\}$, we have
	\begin{align} \label{C:4}
	C_m(n) = \sum_{j=-\infty}^\infty (-1)^j u_m\big(n-j(3j-1)/2\big).
	\end{align}
\end{corollary}

On the other hand, from Theorems \ref{Th2} and \ref{Th3} we easily deduce 
{the following decomposition of $u_m(n)$.}

\begin{corollary}\label{Cor5}
	Let $k$ and $n$ be positive integers. For $m\in\{0,1,2\}$, we have
	\begin{align} \label{C:5b}
		u_m\big(n\big) = \sum_{j=0}^{\infty} C_m\big(n+k(3k+1)/2-j\big)\,\big(M_k(j)+\widetilde{P}_k(j)\big).
	\end{align}
\end{corollary}	
	
In the context given by Corollary \ref{Cor5}, it would be very appealing to have combinatorial
interpretations of
$$\sum_{j=0}^n C_m(j)\,M_k(n-j)$$
and
$$\sum_{j=0}^n C_m(j)\,\widetilde{P}_k(n-j)$$
when $m\in\{0,1,2\}$.

The remainder of the paper is organized  as follows. In Section \ref{S2}, we prove Theorem \ref{Th1}. In Sections \ref{S3} and \ref{S4}, we consider two truncated forms of Euler's pentagonal number theorem in order to prove Theorems \ref{Th2} and \ref{Th3}.  In Section \ref{S5} we prove Corollary \ref{decomp} combinatorially and in Section \ref{S6} we give a combinatorial proof of the Xia-Zhao Truncated Euler pentagonal number theorem used to derive Theorem \ref{Th3}.

\section{Proof of Theorem \ref{Th1}}
\label{S2}

\subsection{Analytic proof}

Here and throughout the paper, we use the following customary $q$-series notation:
\begin{align*}
& (a;q)_n = \begin{cases}
1, & \text{for $n=0$,}\\
(1-a)(1-aq)\cdots(1-aq^{n-1}), &\text{for $n>0$;}
\end{cases}\\
& (a;q)_\infty = \lim_{n\to\infty} (a;q)_n.
\end{align*}
Because the infinite product $(a;q)_{\infty}$ diverges when $a\neq 0$ and $|q| \geqslant 1$, whenever
$(a;q)_{\infty}$ appears in a formula, we shall assume $|q| < 1$.

It is not difficult to see that the number of unimodal compositions of $n$ with largest term $k$ is the coefficient of $q^n$ in $$\frac{q^k}{(q;q)_{k-1}(q;q)_k}.$$
Hence
\begin{align}
\sum_{n=0}^\infty u_0(n)\,q^n = \sum_{k=1}^\infty \frac{q^k}{(q;q)_{k-1}(q;q)_k}.\label{gf0}
\end{align}
More details about this generating function can be found in \cite{Auluck} and \cite[Sec. 2.5]{Stanley}.
Analogously, if $m>0$ and $k\geqslant 0$, the number of unimodal compositions of $n+m$ in which the  largest term is $k+1$ and it occurrs exactly $m$ times  is the coefficient of $q^n$ in
$$\frac{q^{mk}}{(q;q)^2_k}.$$
Thus, for $m=1,2$, we deduce that
\begin{align}
\sum_{n=0}^\infty u_m(n)\,q^n = \sum_{k=0}^\infty \frac{q^{mk}}{(q;q)^2_k}.\label{gf1}
\end{align}
Considering \eqref{gf0} and \eqref{gf1}, we have
\begin{align*}
\sum_{n=0}^\infty \big(u_1(n)-u_2(n)\big)\,q^n = \sum_{k=0}^\infty \frac{q^{k}(1-q^{k})}{(q;q)^2_k} = \sum_{n=0}^\infty u_0(n)\,q^n.
\end{align*}
This concludes the proof.

\subsection{Combinatorial proof}

If $n=0$ the identity of the theorem is true by definition. For the remainder of the proof, we consider $n>0$.  We denoteby $\mathcal U_0(n)$  the set of unimodal compositions of $n$. 
For $m>0$, we denote by $\mathcal U_m(n)$ the set of unimodal compositions of $n+m$ where the maximal part appears exactly $m$ times.

We create an injection $\varphi:\mathcal U_2(n)\to \mathcal U_1(n)$ as follows. If $\lambda\in \mathcal U_2(n)$ with largest parts $\lambda_j =\lambda_{j+1}$, we define $\varphi(\lambda)=\mu$, where $\mu_i=\lambda_i$ if $i\neq j+1$ and $\mu_{j+1}=\lambda_{j+1}-1$.  Since $\lambda$ is a composition of size at least $3$, it follows that  $\lambda_j=\lambda_{j+1}\geqslant 2$. Thus $\varphi(\lambda)\in \mathcal U_1(n)$ and $\mu_{j+1}=\mu_j-1$. The image $\varphi(\mathcal U_2(n))$ of $\varphi$ is the subset of $\mathcal U_1(n)$ consisting of compositions $\mu$ of $n+1$ such that if $m_j$ is the unique largest part, then $\mu_{j+1}=\mu_j-1\geqslant 1$, i.e., in the Ferrers diagram of $\mu$, there is at least one row of length $\mu_j-1$ below the $j$th row.  Clearly the map $\varphi:\mathcal U_2(n) \to \varphi(\mathcal U_2(n))$ is invertible and $\varphi$ is injective. 

Next, we create an bijection $\psi: \mathcal U_1(n)\setminus \varphi(\mathcal U_2(n)) \to \mathcal U_0(n)$. The set $\mathcal U_1(n)\setminus \varphi(\mathcal U_2(n))$ consists of unimodal compositions of $n+1$ with unique largest part $\mu_j$ such that $\mu_i<\mu_j-1$ if $i>j$. We define $\psi(\mu)=\nu$, where $\nu_i=\mu_i$ is $i\neq j$ and $\nu_j=\mu_j-1$. Since $\mu$ is a composition of size at least $2$, we must have $\mu_j>1$ and thus $\nu_j\geqslant 1$ is the largest part of $\nu$. Notice that $\nu$ could have more parts equal to $\nu_j$ but $\nu_j$ is the last part that is largest in $\nu$. 

The inverse of $\psi$ maps $\eta\in \mathcal U_0(n)$ to a composition in  $\mathcal U_1(n)\setminus \varphi(\mathcal U_2(n))$ by adding one to the last of the largest parts of $\eta$. Therefore $\psi$ is a bijection and $u_0(n)=u_1(n)-u_2(n)$. 

\begin{remark} Given a composition of $n$, by adding one to a maximal part, one obtains a composition in $\mathcal U_1(n)$. Thus,  $u_1(n)$ is also equal to the total number of maximal parts in all unimodal compositions of $n$. \end{remark}

\section{Proof of Theorem \ref{Th2}}
\label{S3}


In \cite{Garvan}, as a precursor to the crank of a partition, Garvan defined the crank of a vector partition as follows. 
	\begin{definition} Let $\pi=(\pi_1, \pi_2, \pi_3)$ be a triple of partitions such that $\pi_1$ has distinct parts and $\pi_2$ and $\pi_3$ are unrestricted. Then the crank of $\pi$ is defined as $r(\pi)=\ell(\pi_2)-\ell(\pi_2)$. 
	\end{definition}
	
	The triple $\pi=(\pi_1, \pi_2, \pi_2)$  above is also called a vector partition. It is said to be a vector partition of $n$ (or of  size $n$), denoted $|\pi|=n$, if $|\pi_1|+|\pi_2|+|\pi_3|=n$. 
	
	Then the weighted number of vector partitions of $n$ and crank $k$ is $$N_V(k,n)=\sum_{|\pi|=n, r(\pi)=k}(-1)^{\ell(\pi_1)}.$$
	
	In \cite{Andrews88} it is shown that, if $n\neq 1$, $N_V(k,n)$ equals the number of ordinary partitions of $n$ with crank $k$. 

From \cite[Theorem 7.19]{Garvan}, for a given integer $k$, the generating function for $N_V(k,n)$ 
is given by

\begin{align}
\frac{1}{(q;q)_\infty} \sum_{n=1}^\infty (-1)^{n-1} q^{n(n-1)/2+n|k|}(1-q^n). \label{gfc}
\end{align}
Summing \eqref{gfc} over $k>0$ yields
\begin{align*}
\sum_{n=0}^\infty C_0(n)\,q^n = \frac{-1}{(q;q)_\infty} \sum_{n=1}^\infty (-1)^{n} q^{n(n+1)/2}.
\end{align*}
In a similar way, summing \eqref{gfc} over $k\geqslant 0$ yields
\begin{align*}
\sum_{n=0}^\infty C_1(n)\,q^n = \frac{1}{(q;q)_\infty} \sum_{n=0}^\infty (-1)^n q^{n(n+1)/2}.
\end{align*}

In \cite[Sec 2.5]{Stanley}, Stanley provided the following identities:
\begin{align*}
\sum_{n=0}^\infty u_0(n)\,q^n = \frac{-1}{(q;q)^2_\infty} \sum_{n=1}^\infty (-1)^{n} q^{n(n+1)/2}
\end{align*}
and
\begin{align*}
\sum_{n=0}^\infty u_1(n)\,q^n = \frac{1}{(q;q)^2_\infty} \sum_{n=0}^\infty (-1)^n q^{n(n+1)/2}.
\end{align*}
For $m\in\{0,1,2\}$, we easily deduce that
$$
\sum_{n=0}^\infty u_m(n)\,q^n = \frac{1}{(q;q)_\infty} \sum_{n=0}^\infty C_m(n)\,q^n.
$$

In \cite{Andrews12}, Andrews and Merca considered Euler's pentagonal number theorem 
$$
(q;q)_\infty = \sum_{n=-\infty}^\infty (-1)^n q^{n(3n-1)/2},
$$	
and they proved the following truncated form: 

\begin{theorem}[Andrews-Merca] \label{AM} for any $k\geqslant 1$, 
\begin{equation} \label{TPNT}
\frac{(-1)^{k-1}}{(q;q)_\infty} \sum_{n=1-k}^{k} (-1)^{n} q^{n(3n-1)/2}= (-1)^{k-1}+ \sum_{n=k}^\infty \frac{q^{{k\choose 2}+(k+1)n}}{(q;q)_n}
\begin{bmatrix}
n-1\\k-1
\end{bmatrix},
\end{equation}
where  
$$
\begin{bmatrix}
n\\k
\end{bmatrix} 
=
\begin{cases}
\dfrac{(q;q)_n}{(q;q)_k(q;q)_{n-k}}, &  \text{if $0\leqslant k\leqslant n$},\\
0, &\text{otherwise.}
\end{cases}
$$\end{theorem}
We note that the series on the right hand side  of \eqref{TPNT} is the generating function for $M_k(n)$, i.e.,
\begin{equation}
\sum_{n=0}^\infty M_k(n) q^n = \sum_{n=k}^\infty \frac{q^{{k\choose 2}+(k+1)n}}{(q;q)_n}
\begin{bmatrix}
n-1\\k-1
\end{bmatrix}. \label{mk}
\end{equation}

For $m\in\{0,1,2\}$, multiplying both sides of \eqref{TPNT} by
$$
\sum_{n=0}^\infty C_m(n)\,q^n 
$$
we obtain
\begin{align*}
(-1)^{k-1}\left( \Big(\sum_{n=0}^\infty u_m(n)\, q^n \Big) \Big( \sum_{n=1-k}^k (-1)^n q^{n(3n-1)/2} 
\Big) - \sum_{n=0}^\infty C_m(n)\, q^{n}  \right)  \\
= \left( \sum_{n=0}^\infty C_m(n)\, q^{n} \right) \left( \sum_{n=0}^\infty M_k(n)\, q^n\right).
\end{align*}
The proof of Theorem \ref{Th2} follows easily considering Cauchy's multiplication of two power series.

\section{Proof of Theorem \ref{Th3}}
\label{S4}

The proof of this identity is quite similar to the proof of Theorem \ref{Th2}.
In \cite{Xia},  Xia and Zhao considered Euler's pentagonal number theorem 
and proved the following truncated form:
\begin{equation} \label{TPNT2}
\frac{1}{(q;q)_\infty} \sum_{n=-k}^{k} (-1)^{n}\, q^{n(3n-1)/2}= 1+(-1)^k\frac{q^{k(k+1)/2}}{(q;q)_k} \sum_{n=0}^\infty \frac{q^{(n+k+1)(k+1)}}{(q^{n+k+1};q)_\infty}.
\end{equation}

Multiplying both sides of \eqref{TPNT2} by 
$$
\sum_{n=0}^\infty C_m(n)\,q^n,
$$
we obtain
\begin{align*}
& (-1)^{k} \left( \bigg( \sum_{n=1}^\infty u_m(n)\, q^n \bigg) \bigg(  \sum_{n=-k}^{k} (-1)^{n}\, q^{n(3n-1)/2}\bigg) -\sum_{n=0}^\infty C_m(n)\,q^n\right)   \\
& =   \left( \sum_{n=0}^\infty C_m(n)\,q^n \right) \left( \sum_{n=0}^\infty \widetilde{P}_k(n)\,q^n\right),
\end{align*}
where we have invoked 
the generating function for $\widetilde{P}_k(n)$ given in \cite{Xia},
$$
\sum_{n=0}^\infty \widetilde{P}_k(n)\, q^n = \frac{q^{k(k+1)/2}}{(q;q)_k} \sum_{n=0}^\infty \frac{q^{(n+k+1)(k+1)}}{(q^{n+k+1};q)_\infty}.
$$
The proof follows easily considering Cauchy's multiplication of two power series.


\section{Combinatorial proof of Corollary \ref{decomp}} \label{S5}

We use the fact that, for $n\neq 1$, $N_V(k,n)$, the weighted number of vector partitions of $n$ and crank $k$,   equals the number of ordinary partitions of $n$ with crank $k$.
A combinatorial proof of this fact was given by Dyson in \cite{Dyson}. 

For $n=1$, the righthand side of \eqref{C:4} equals $u_m(1)-u_m(0)$. Since the only unimodal composition of $1$ is $(1)$, we have $u_0(1)=u_1(1)=1$ and $u_2(1)=0$ and the statement can be verified directly.  

If $n>1$, we use vector partitions and $N_V(n,k)$ to count partitions of $n$ with crank $k$. 

First, we simpify the definition of $N_V(n,k)$. For every positive integer $j$, we denote the corresponding pentagonal partitions by $G_j=(2j-1, 2j-2, \ldots, j+1, j)$ and $G_{-j}= (2j, 2j-1, \ldots, j+2, j+1))$. Moreover,  $G_0=\emptyset$. For $j \in \mathbb Z$, we have  $|G_j|=j(3j-1)/2$.
Let $\mathcal V_e(n)$, respectively $\mathcal V_o(n)$, be  the set of vector partitions $\pi=(\pi_1, \pi_2,\pi_3)$ with $\ell(\pi_1)$ even, respectively odd,  and $\pi_1$ is not a pentagonal partition. 
 We create a bijection from $\mathcal V_o(n)$ to $\mathcal V_e(n)$ by mapping $(\pi_1, \pi_2, \pi_3)$ to $(\varphi(\pi_1), \pi_2, \pi_3)$, where $\varphi$ is 
 Franklin's involution for the combinatorial proof of Euler's pentagonal number theorem. Then, $$N_V(n,k)= \sum_{j=-\infty}^\infty (-1)^j|\{\pi=(G_j, \pi_2, \pi_3) \mid \ell(\pi_2)-\ell(\pi_3)=k\}|.$$
 
 Before we continue with the proof of the corollary we introduce more notation. Given a partition $\lambda$, we denote by $\lambda'$ the conjugate partition of $\lambda$ and by $\lambda^*$ the composition whose parts are exactly the parts of $\lambda'$ written in non-decreasing order. Note that $\lambda^*$ is the composition whose Ferrers diagram is obtained by rotating the Ferrers diagram of $\lambda$ by $90^\circ$ counterclockwise. For example, if $\lambda=(5,4,4,1,1,1)$, then $\lambda'=(6,3,3,3,1)$ and $\lambda^*=(1,3,3,3,6)$. 
 
 First, we prove \eqref{C:4} for $m=2$. For each $j \in \mathbb Z$, we create a bijection from $\{\pi=(G_j, \pi_2, \pi_3) \mid |\pi|=n, \ell(\pi_2)=\ell(\pi_3)\}$ to $\mathcal U_2(n-j(3j-1)/2)$, the set of unimodal compositions of $n+2-j(3j-1)/2$ with exactly two maximal parts, as follows. Start with $(G_j, \pi_2, \pi_3)$ with $\ell(\pi_2)=\ell(\pi_3)$ and  add a part equal to $1$ to each of $\pi_2$ and $\pi_3$ to obtain partitions $\alpha_2=\pi_2\cup\{1\}$ and $\alpha_3=\pi_3\cup\{1\}$, respectively. Define a unimodal composition $\psi(\pi)\in \mathcal U_2(n-j(3j-1)/2)$ by the sequence whose initial parts are the parts of $\alpha_3^*$ followed by the parts of $\alpha_2'$.  For the inverse, start with a unimodal composition $\eta\in \mathcal U_2(n-j(3j-1)/2)$ and suppose the maximal parts of $\eta$ are $\eta_i$ and $\eta_{i+1}$.  Notice that $\eta_{i-1}=0$ or $\eta_{i-1}<\eta_i$ and $\eta_{i+2}=0$ or $\eta_{i+2}<\eta_{i+1}$. Subtract one from each of $\eta_i$ and $\eta_{i+1}$. Then, parts $(\eta_1, \eta_2, \ldots \eta_{i-1}, \eta_i-1)$ form  a composition which  is $\pi_3^*$ for a partition $\pi_3$ with $\ell(\pi_3)=\eta_i-1$. Moreover, parts $(\eta_{i+1}-1, \eta_{i+2}, \ldots, \eta_{\ell(\eta)})$ form  a partition which is the conjugate of partition $\pi_2$ with $\ell(\pi_2)=\eta_{i+1}-1=\eta_i-1$. Thus $\ell(\pi_2)=\ell(\pi_3)$,  $\pi=(G_j, \pi_2, \pi_3)$ has size $n$, and $\psi(\pi)=\eta$. 
 
 Therefore  \begin{align*} C_2(n) & =\sum_{j=-\infty}^\infty (-1)^j|\{\pi=(G_j, \pi_2, \pi_3) \mid \ell(\pi_2)=\ell(\pi_3)\}|\\ & =  \sum_{j=-\infty}^\infty (-1)^j u_2\big(n-j(3j-1)/2\big).\end{align*}

Next, we prove \eqref{C:4} for $m=0$. We have $$C_0(n)=\sum_{j=-\infty}^\infty (-1)^j|\{\pi=(G_j, \pi_2, \pi_3) \mid \ell(\pi_2)>\ell(\pi_3)\}|.$$ As above, for each $j \in \mathbb Z$,  we create a bijection from $\{\pi=(G_j, \pi_2, \pi_3) \mid |\pi|=n, \ell(\pi_2)>\ell(\pi_3)\}$ to $\mathcal U_0(n-j(3j-1)/2)$ as follows.  Start with  $\pi=(G_j, \pi_2, \pi_3)$, where $\ell(\pi_2)>\ell(\pi_3)$    and define $\psi(\pi)$ to be the composition whose initial parts are the parts of $\pi_3^*$ followed by the parts of $\pi_2'$.  The first part of $\pi_2'$ is the first maximal part of $\psi(\pi)$. To see that the transformation is invertible, start with a unimodal composition $\eta$ of $n-j(3j-1)/2$ and suppose $\eta_i$ is the first maximal part. Then the first $i-1$ parts of $\eta$ form $\pi_3^*$ for a partition $\pi_3$, and the remaining parts (in the given order) form $\pi_2'$ for a partition $\pi_2$. Then $\pi=(G_j, \pi_2, \pi_3)$ has size $n$, $\ell(\pi_2)>\ell(\pi_3)$, and $\psi(\pi)=\eta$. 

Finally, we prove \eqref{C:4} for $m=1$. We have $$C_1(n)=\sum_{j=-\infty}^\infty (-1)^j|\{\pi=(G_j, \pi_2, \pi_3) \mid \ell(\pi_2)\geqslant\ell(\pi_3)\}|.$$ As above, for each $j \in \mathbb Z$,  we create a bijection from $\{\pi=(G_j, \pi_2, \pi_3) \mid |\pi|=n, \ell(\pi_2)\geqslant\ell(\pi_3)\}$ to $\mathcal U_1(n-j(3j-1)/2)$ the set of unimodal compositions of $n+1-j(3j-1)/2$ with exactly one maximal part, as follows. Start with $(G_j, \pi_2, \pi_3)$ with $\ell(\pi_2)\geqslant \ell(\pi_3)$  and  add a part equal to $1$ to $\pi_2$ to obtain a partition $\alpha_2=\pi_2\cup\{1\}$. Define a unimodal composition $\psi(\pi)$ by the sequence whose initial parts are the parts of $\pi_3^*$ followed by the parts of $\alpha_2'$. As in the case $m=2$, this transformation is invertible.

\section{Combinatorial Proof of the Xia-Zhao Truncated Pentagonal Number Theorem} \label{S6}

Recently, Xia and Zhao introduced a new truncated pentagonal number theorem which was used in section \ref{S4} to prove Theorem \ref{Th3}. In this section we give a combinatorial proof of this theorem. 

We denote by $\mathcal P(n)$ the set of partitions of $n$ and by $\widetilde{\mathcal{P}}_k(n)$ the set of partitions in which every integer less than or equal to $k$ appears as a part at least once and the first part larger than $k$ appears at least $k+1$ times. Then $P_k(n)=|\widetilde{\mathcal{P}}_k(n)|$. 

\begin{theorem}[Xia-Zhao] \label{XZ} Let $n$ and $k$ be positive integers. Then, \begin{equation}\label{xz}(-1)^k\sum_{j=-k}^k(-1)^jp(n-j(3j+1)/2)=\widetilde P_k(n).\end{equation} \end{theorem}

We note that there is a slight error in the definition of $\widetilde P_k(n)$ in \cite{Xia} where $k$ is specified to be the least integer satisfying the definition. It is clear from the generating function  used in \cite{Xia} that $k$ does not have to be the smallest integer with the properties in the definition of $\widetilde P_k(n)$.

\begin{proof} First note that for $n=1$ each side of \eqref{xz} is zero. For the remainder of the proof $n>1$. The statement of the theorem is equivalent to \begin{align} \label{k=1}p(n-1)+p(n-2)-p(n)&=\widetilde P_1(n)   \\ \label{genk}p\left(n-\frac{k(3k-1)}{2}\right)+p\left(n-\frac{k(3k+1)}{2}\right)& = \widetilde P_{k-1}(n)+\widetilde P_k(n), \ k\geqslant 2.
\end{align}

To prove \eqref{k=1}, notice that by adding a part equal to $1$ to a partition in $\mathcal P(n-1)$ we obtain a bijection between $\mathcal P(n-1)$ and the set of partitions in $\mathcal P(n)$ that have at least a part equal to $1$. Thus, $p(n)-p(n-1)$ is the number of partitions of $n$ that do not contain a part equal to $1$. Similarly, by adding a part equal to $2$ to a partition in $\mathcal P(n-2)$, we see that $p(n-2)$ is also the number of partitions in $\mathcal P(n)$ that have at least one part equal to $2$. Then, $p(n-1)+p(n-2)-p(n)=p_{1,2}(n)-p_{\overline{1,2}}(n)$, where $p_{1,2}(n)$ is the number of partitions of $n$ that have at least one part equal to each $1$ and $2$ and $p_{\overline{1,2}}(n)$ is the number of partitions of $n$ that do not have parts equal to $1$ or $2$.

Next, for a partition $\lambda$ of $n$ that does not have parts equal to $1$ or $2$, if the smallest part is equal to $a\geqslant 3$, replace the part $a$ by one part equal to $1$ and $(a-1)/2$ parts equal to $2$ if $a$ is odd, and by two parts equal to $1$ and $(a-2)/2$ parts equal to $2$ if $a$ is even. Then, $p(n-1)+p(n-2)-p(n)$ equals $|\mathcal P^*_{1,2}(n)|$, where $\mathcal P^*_{1,2}(n)$ is the set of partitions $\mu$ of $n$ with at least one part equal to $2$ and satisfying 
\begin{itemize}
\item[(i)]$\mu$ has  at least three parts equal to $1$, or
\item[(ii)] $\mu$ has one or two parts equal to $1$,  at least one part greater than $2$, and the sum of all parts equal to $1$ or $2$ is larger than the smallest part greater than $2$. 
\end{itemize}
We create a bijection $f$ between $\mathcal P^*_{1,2}(n)$ and $\widetilde{\mathcal{P}}_1(n)$ as follows. 

If $\mu\in \mathcal P^*_{1,2}(n)$  satisfies (i), we replace two parts equal to $1$ by a part equal to $2$ to obtain a partition in $\widetilde{\mathcal{P}}_1(n)$ that contains (at least two) parts equal to $2$. For example, if $\mu=(8,8,5,3,2,2,1,1,1,1)$, then $f(\mu)=(8,8,5,3,2,2,2,1,1)$. 

If $\mu\in \mathcal P^*_{1,2}(n)$  satisfies (ii) and $b$ is the smallest part greater than $2$, we replace all parts equal to $1$ and $2$ by a part equal to $b$ and parts equal to $1$. We obtain a partition in $\widetilde{\mathcal{P}}_1(n)$ with no part equal to $2$.  For example, if  $\mu=(8,8,5,5, 2,2,2,1,1)$, then $f(\mu) =(8,8,5,5,5,1,1,1)$.

To see that this transformation is invertible, start with $\lambda \in  \widetilde{\mathcal{P}}_1(n)$. If $\lambda$ has  parts equal to $2$, then it has at least two parts equal to $2$. Replace one part equal to $2$ by two parts equal to $1$ to obtain a partition $\mu\in \mathcal P^*_{1,2}(n)$ satisfying (i) such that $f(\mu)=\lambda$. If $\lambda$ does not have parts equal to $2$, suppose $b$ is the smallest part greater than $2$. Then $\lambda$ has at least two parts equal to $b$ and $m$ parts equal to $1$, for some  $m\geqslant 1$. Remove one part equal to $b$ and all $m$ parts equal to $1$. We have $m+b\geqslant 4$. If $m+b$ is even, insert two parts equal to $1$ and $(m+b-2)/2$ parts equal to $2$.  If $m+b$ is odd, insert one part equal to $1$ and $(m+b-1)/2$ parts equal to $2$. We obtain a partition $\mu\in \mathcal P^*_{1,2}(n)$ satisfying (i) such that $f(\mu)=\lambda$.

Let $k\geqslant 2$. To prove \eqref{genk}, we create a bijection $$g:\mathcal P\left(n-\frac{k(3k-1)}{2}\right)\cup\mathcal P\left(n-\frac{k(3k-1)}{2}\right) \to\widetilde{\mathcal{P}}_{k-1}(n)\sqcup  \widetilde{\mathcal{P}}_k(n).$$

Note that $\widetilde{\mathcal{P}}_{k-1}(n)\cap  \widetilde{\mathcal{P}}_k(n)$ consists of partitions of $n$ in which parts $1, 2, \ldots, k-1$ occur at least once, part $k$ occurs at least $k$ times, and the first part greater than $k$ appears at least $k+1$ times. Thus, these partitions occur twice in $\widetilde{\mathcal{P}}_{k-1}(n)\sqcup  \widetilde{\mathcal{P}}_k(n)$.

First, let $\lambda \in \mathcal P\left(n-\frac{k(3k-1)}{2}\right)$. We define $g(\lambda)$ to be the partition of $n$ obtained  from $\lambda$ by inserting  one part equal to each $i$ for $1\leqslant i\leqslant k-1$ and  $k$ parts equal to $k$. Then, $g(\lambda)$ is a partition in $\widetilde{\mathcal{P}}_{k-1}(n)$ in which the first part larger than $k-1$ is $k$.  Note that if the first part greater than $k$ in $g(\lambda)$ occurs at least $k+1$ times, then $g(\lambda)$  is also a partition in $\widetilde{\mathcal{P}}_{k}(n)$.

Next, let $\lambda \in \mathcal P\left(n-\frac{k(3k+1)}{2}\right)$ and denote by $m_k\geqslant 0$ be the multiplicity of $k$ in $\lambda$. 

\underline{Case 1:} $\lambda$ has no part greater than $k$ or $x$ is  the first part of $\lambda$ greater than $k$  and $x\geqslant k+m_k+1$. Note that $m_k$ could be $0$.  We define $g(\lambda)$ to be the partition of $n$ obtained from $\lambda$ by removing all $m_k$ parts equal to $k$ and inserting one part equal to each $i$ for $1\leqslant i\leqslant k-1$ and also inserting $k$ parts equal to $k+m_k+1$. Then, $g(\lambda)$ is a partition in $\widetilde{\mathcal{P}}_{k-1}(n)$ in which the first part larger than $k-1$ is equal to $k+m_k+1$. In particular, $g(\lambda)$ has no part equal to $k$. 

\underline{Case 2:}  $\lambda$ has a part greater than $k$ and the first part greater than $k$ is $x\leqslant k+m_k$. Then, we must have $m_k>0$.  We define $g(\lambda)$ to be the partition of $n$ obtained from $\lambda$ by removing all $m_k$ parts equal to $k$ and
 inserting one part equal to each $i$ for $1\leqslant i\leqslant k-1$ and  also inserting $k+m_k+1-x$ parts equal to $k$ and  $k$ parts equal to $x$. Then, $g(\lambda)$ is a partition in $\widetilde{\mathcal{P}}_{k}(n)$. Note that, if $m_k\geqslant x-1$, then $k$ appears in $g(\lambda)$ at least $k$ times and $g(\lambda)$ is also a partition in $\widetilde{\mathcal{P}}_{k-1}(n)$.

Thus we defined  a mapping $g$ and elements of $\widetilde{\mathcal{P}}_{k-1}(n)\cap  \widetilde{\mathcal{P}}_k(n)$ appear (at least) twice in the image of $g$. 

Next we show that the transformation $g$, as a mapping on the disjoint union $\widetilde{\mathcal{P}}_{k-1}(n)\sqcup  \widetilde{\mathcal{P}}_k(n)$, is invertible. Let $\mu \in \widetilde{\mathcal{P}}_{k-1}(n)\cup  \widetilde{\mathcal{P}}_k(n)$. We denote by $\ell_k$ the multiplicity of $k$ in $\mu$. 

\underline{Case I:} $\mu \in \widetilde{\mathcal{P}}_k(n)$. We denote by $y$ the first part of $\mu$ larger than $k$. Thus, $y$ occurs in $\mu$ at least $k+1$ times and $\ell_k\geqslant 1$.  In $\mu$, remove one part equal to each $i$ for $1\leqslant i\leqslant k-1$ and also remove $k$ parts equal to $y$, and inset $y-k-1$ parts equal to $k$. We obtain a partition $\lambda\in  \mathcal P\left(n-\frac{k(3k+1)}{2}\right)$ from Case 2 above, i.e., $y$ is the first part of $\lambda$  greater than $k$ and $y\leqslant k+m_k$, where $m_k\geqslant 1$ is the multiplicity of $k$ in $\lambda$. Then, $g(\lambda)=\mu$. 

\underline{Case II(a):} $\mu \in \widetilde{\mathcal{P}}_{k-1}(n)$ and $\mu$ has parts equal to $k$. Thus, $\mu$ has at least $k$ parts equal to $k$. We remove from $\mu$ one part equal to each $i$ for $1\leqslant i\leqslant k-1$ and we also remove $k$ parts equal to $k$. We obtain a partition $\lambda \in  \mathcal P\left(n-\frac{k(3k-1)}{2}\right)$ such that $g(\lambda)=\mu$. 

\underline{Case II(b):} $\mu \in \widetilde{\mathcal{P}}_{k-1}(n)$ and $\mu$ has no parts equal to $k$.  We denote by $y$ the first part of $\mu$ larger than $k$. Thus, $y$ occurs in $\mu$ at least $k$ times. In $\mu$, remove one part equal to each $i$ for $1\leqslant i\leqslant k-1$ and also remove $k$ parts equal to $y$, and insert $y-k-1$ parts equal to $k$. We obtain a partition $\lambda\in  \mathcal P\left(n-\frac{k(3k+1)}{2}\right)$ from Case 1 above, i.e., there is no part greater than $k$ or $y$ is the first part of $\lambda$  greater than $k$ and $y\geqslant k+m_k+1$, where $m_k\geqslant 1$ is the multiplicity of $k$ in $\lambda$. Then, $g(\lambda)=\mu$. 

Note that for $\lambda \in \widetilde{\mathcal{P}}_{k-1}(n)\cap  \widetilde{\mathcal{P}}_k(n)$, we transform one copy of $\lambda$ according to the rule in Case I to obtain a partition of $n-k(3k+1)/2$ and we transform another  copy of $\lambda$ according to the rule in Case II(a) to obtain a partition of $n-k(3k-1)/2$. 
\end{proof}

\begin{remark} Combinatorially, Theorem \ref{AM} means that, for $k, n \geqslant 1$, we have \begin{equation}\label{am}(-1)^k\sum_{j=1-k}^k(-1)^jp(n-j(3j+1)/2)=M_k(n).\end{equation}
\end{remark} Then, \eqref{am} and \eqref{xz} imply that, for $k, n\geqslant 1$, we have  \begin{equation}\label{MP} p(n-k(3k+1)/2)=M_k(n)+  \widetilde{\mathcal{P}}_k(n).\end{equation}

If we denote by $p(N,M,n)$ the number of partitions of $n$ with at most $M$ parts, each at most $N$, then \cite[Theorem 3.1]{Andrews88} gives a combinatorial proof for the generating function for $p(N,M,n)$, i.e., \begin{equation}\label{qbin}\sum_{n=0}^\infty p(N,M,n) q^n= \frac{(q;q)_{N+M}}{(q;q)_N(q;q)_M}.\end{equation} Using this fact, one can easily show  combinatorially that $M_k(n)$ equals the number of partitions in $\mathcal P_{k-1}(n)$ that do not have parts equal to $k$. Then, the combinatorial proof of Theorem \ref{XZ} gives a new combinatorial proof of Theorem \ref{AM}. However, the combinatorial proof of \eqref{qbin} is not bijective. Rather, in \cite{Andrews88}, it is proved that $p(N,M,n)$ and the  sequence whose generating function is the right hand side of \eqref{qbin} satisfy the same recurrence and have the same initial values. Thus, our combinatorial proof of Theorem \ref{AM} is not bijective and neither are the known proofs so far. A bijective proof that $M_k(n)$ equals the number of partitions in $\mathcal P_{k-1}(n)$ that do not have parts equal to $k$, would lead to  a  bijective proof of Theorem \ref{AM}. This remains an open problem.

\bigskip


\end{document}